\documentclass[11pt]{article}
\usepackage{amsthm,amsfonts,amsmath,amssymb}
\usepackage{multirow}
\newtheorem{thm}{Theorem}[section]
\newtheorem{prop}[thm]{Proposition}
\newtheorem{lem}[thm]{Lemma}

\newtheorem{rem}[thm]{Remark}

\def\a{\alpha} \def\b{\beta}   
    
 \def\ld{\lambda}  
\def\D{\Delta}

 \def\lg{\langle} \def\rg{\rangle}
\def\nd{\mathrel{\bigm|\kern-.7em/}}

\def\f{\noindent}

\def\Aut{\hbox{\rm Aut\,}}

\def\Syl{\hbox{\rm Syl}}

\def\mod{\hbox{\rm mod }}

\def\demo{\f{\bf Proof}\hskip10pt}

\def\qed{\hfill $\Box$}

\def\Q{{\rm Q}}
\def\D{{\rm D}}
\def\SD{{\rm SD}}
\def\C{{\rm C}}

\def\lg{\langle}
\def\rg{\rangle}

\begin{document}
\title{Finite groups which have maximal covers
\thanks{This work was supported by NSFC (no. 11971280, 11771258). }}
\author{Lifang Wang and Lijian An\thanks{corresponding author. e-mail: anlj@sxnu.edu.cn}\\
Department of Mathematics, Shanxi Normal University\\
Linfen, Shanxi, 041004 P. R. China}
\date{}
 \maketitle
\begin{abstract}
  Let $\ld(G)$ be the maximum number of subgroups in an irredundant covering of
  a finite group $G$. We prove that the finite groups with
  $\ld(G)=|G|-t$, where $t\leq 5$, are solvable, and classify such groups.

\medskip


\end{abstract}

\baselineskip=16pt

\section{Introduction}

\ \ \ Let $G$ be a finite group. A {\it cover} of $G$ is a collection
of subgroups of $G$ whose union is $G$. The cover is {\it irredundant}
if no proper sub-collection is also a cover.
Many different aspects of covers have been studied.
One of these aspects is to study
the structure of groups by using the number of subgroups in a cover of a group.
Cohn \cite{Cohn}
defined $\sigma(G)$ to be the smallest integer $n$ such that the group $G$ is
covered by $n$ proper subgroups. He proved some properties of $\sigma(G)$
and described the structure of the groups with
small $\sigma(G)$.

Since then, many scholars study groups by $\sigma(G)$. On one hand,
some scholars calculated $\sigma(G)$ for given groups.
For example, Cohn \cite{Cohn} proved
$\sigma(A_5)=10$ and $\sigma(S_5)=16$.
Tomkinson \cite{Tomkinson1997} proved that
$\sigma(G)=1+p^a$, where $p^a$ is the order of a particular chief factor of $G$, if $G$ is a finite solvable group.
 Holmes \cite{P.E.Holmes2006} determined $\sigma(G)$
for $G=M_{11}, M_{22}$ and $ M_{23}$ respectively. Kappe et al. \cite{KNS2016}
determined $\sigma(S_n)$ for $n=8, 9, 10, 12$.
On the other hand, some scholars determined which numbers can occur as
 $\sigma(G)$
for groups, and, when possible, to characterize the
groups having the same value of $\sigma(G)$.
For example,
Cohn \cite{Cohn} described the structure of the groups with
$\sigma(G)=3,4,5$. Abdollahi et al. \cite{AAJM2005} gave a description of finite groups with $\sigma(G)=6$.
Tomkinson \cite{Tomkinson1997} proved that
there is no finite groups $G$ such that $\sigma(G)=7$.
Zhang \cite{Zhang2006} proved that there is no finite
groups $G$ satisfying $\sigma(G)=11$ or $\sigma(G)=13$.
For more results, the readers refer to \cite{AJ2007, BFM2015, Serena2003, Swartz2016}.

As the dual of $\sigma(G)$, Rog${\rm\acute{e}}$rio
\cite{Rogerio} defined $\ld(G)$  as the maximum number of subgroups
in an irredundant covering of a group $G$.
In \cite{Rogerio}, the author gave some basic properties of $\ld(G)$,
calculated $\ld(G)$ for some groups and
classified the finite groups with $\sigma(G)=\ld(G)$ and $\ld(G)=3,4,5$, respectively.
Bastos et al. \cite{BLR2020} described the structure of groups $G$
with $\ld(G)=6$.

We continue to investigate the structure of finite groups
by $\ld(G)$. In \cite{BLR2020,Rogerio}, the authors studied the groups with small $\ld(G)$.
In the paper, we investigate the groups with large $\ld(G)$ and classify the groups with $\ld(G)=|G|-t$, where $t\leq 5$.

The notation and terminology are standard, see \cite{Ber1}.
 We
use ${\rm C}_n, {\rm D}_{2^n}, {\rm Q}_{2^n}$ and ${\rm C}_p^n$ to denote the
cyclic group of order $n$, the dihedral group of order $2^n$, the
generalized quaternion group of order $2^n$ and the elementary abelian group of order $p^n$,
respectively.

Let $G$ be a finite $p$-group.
For any positive integer $s$, we define
$$\Omega_s(G)=\lg a\in G\mid a^{p^s}=1\rg\ {\rm and}\ \mho_s(G)=\lg
a^{p^s}\mid a\in G\rg.$$
$$\Omega_{\{s\}}(G)=\{ a\in G\mid a^{p^s}=1\}\ {\rm and}\
\mho_{\{s\}}(G)=\{ a^{p^s}\mid a\in G\}.$$

If $a,b$ are two elements of $G$,  the commutator of $a,b$
is defined as $[a,b]=a^{-1}b^{-1}ab$.
If $H$ and $K$ are subgroups of $G$ with $G=HK$ and $[H,K]=1$, we call $G$ a
central product of $H$ and $K$, denoted by $G=H\ast K$. Clearly, $H\cap K\leq Z(G)$.

\section{Preliminaries}

\ \ \ In this section, we list some results which often are used in this paper.
The following Lemma tells us $\ld(G)$ is the number of maximal cyclic subgroups of a finite group $G$.

\begin{lem}{\rm \cite[Proposition 4]{Rogerio}}\label{Lemma=ld(G)=the number of maximal cyclic sub}
Let $G$ be a finite group and $\langle x_i\rangle, i=1,2,\cdots,n$ be the maximal cyclic
subgroups of $G$. Then $\bigcup_{i=1}^{n}\langle x_i\rangle $ is an irredundant covering of $G$ and $\lambda(G)=n$.
\end{lem}

\begin{lem}{\rm \cite[Proposition 5]{Rogerio}}\label{Lemma=property}
Let $G$ be a finite group. \\
{\rm (i)} $\ld(H)\leq\ld(G)$ for any subgroup $H$ of $G$.\\
{\rm (ii)} If $N\unlhd G$, then $\ld(G/N)\leq\ld(G)$. For
$N=\cap_{i=1}^{\ld}H_i$, $\ld(G/N)=\ld(G),$
where $\ld=\ld(G)$ and $H_1, H_2,\cdots,H_{\ld}$ are the
maximal cyclic subgroups of $G$.

\end{lem}

\begin{lem}{\rm \cite[Proposition 6 and Proposition 7]{Rogerio}}\label{Lemma=ld(C_p times C_p^k)ld(D_2n)ld(Q_2^n)}\\
$(1)\ \ld(\C_p\times\C_{p^k})=kp-k+2;$\\
$(2) \ \ld(\D_{2n})=n+1;$\\
$(3) \ \ld(\Q_{2^n})=2^{n-2}+1.$
\end{lem}

\begin{lem}\label{Lemma=ld(D_8times C_2)}
$\ld(\D_8\times \C_2)=12.$
\end{lem}

\demo Let $G=\D_8\times \C_2$. Assume without loss generality
$$G=\lg a,b,c\mid a^4=b^2=c^2=1, [a,b]=a^2, [b,c]=[a,c]=1\rangle.$$
It is clear that $\exp(G)=4$. Hence, the order of maximal cyclic subgroups is $2$ or $4$.
By Lemma \ref{lem=commutator gongshi}, we have
$$(a^ib^jc^k)^2=(a^ib^j)^2=a^{2i}[a^i, b^{-j}]=a^{2i(1-j)}. \eqno(*)$$
Hence,
$$(a^ib^jc^k)^2=1\Longleftrightarrow a^{2i(1-j)}=1
\Longleftrightarrow i\equiv 0(\mod 2) \ \ {\mbox{or}}\ \ j\equiv 1(\mod 2).$$
It follows that
$$\Omega_{\{1\}}(G)=\{1,b,c, bc, a^2, a^2b, a^2c, a^2bc, ab, a^3b, abc, a^3bc\}.$$
Thus, the number of maximal cyclic subgroups of order $4$ is
$\frac{|G|-|\Omega_{\{1\}}(G)|}{2}=\frac{16-12}{2}=2$.

By $(*)$, $\mho_{\{1\}}(G)=\{1,a^2\}$. Hence, the
number of maximal cyclic subgroups of order $2$ is
$|\Omega_{\{1\}}(G)|-|\mho_{\{1\}}(G)|=12-2=10$.
Therefore, $\ld(G)=10+2=12$. \qed

\begin{lem}\label{lem=commutator gongshi}{\rm (\cite{Xu})}
Assume $G$ is a metabelian $p$-group. Let $a,b$ be elements of
$G$ and $n$ a positive integer. Then
$$(ab^{-1})^n=a^n\prod_{i+j\le n} [ia,jb]^{n\choose {i+j}}b^{-n},$$
where $[ia,jb]=[a,b,\underbrace{a,\ldots ,a}_{i-1},
\underbrace{b,\ldots ,b}_{j-1}].$
\end{lem}

\section{Some properties}

\ \ \ By Lemma \ref{Lemma=ld(G)=the number of maximal cyclic sub},
$\ld(G)$ is the number of maximal cyclic subgroups of $G$.
In the following, we assume that
$$G=\lg g_1\rg\cup\lg g_2\rg\cup\cdots\cup\lg g_{\ld(G)}\rg
\ \ {\mbox{\rm and}}\ \ o(g_1)\geq o(g_2)\geq\cdots\geq o(g_{\lambda(G)}),$$
where $\langle g_1\rangle, \lg g_2\rg,\cdots, \lg g_{\lambda(G)}\rg$ are all maximal cyclic subgroups of $G$.

\begin{prop} \label{Prop=o(g_1)=3}
{\rm (1)}  $o(g_1)=2$ if and only if $\ld(G)=|G|-1$.

{\rm (2)} If $o(g_1)=3$, then $\ld(G)=|G|-k-1$, where $k$ is the number of subgroups of order $3$ of $G$.
\end{prop}

\demo (1) It is clear that $o(g_1)=2$ if and only if $G\cong C_2^n$.
The result follows.

(2) Since $o(g_1)=3$,
each subgroup of order $3$ of $G$ is maximal cyclic in $G$. Assume that $\lg g_1\rg, \lg g_2\rg,\cdots \lg g_k\rg$
are all subgroups of order $3$ of $G$.
Let $$G=\lg g_1\rg\cup \cdots \lg g_k\rg\cup\lg g_{k+1}\rg\cdots\lg g_{\ld(G)}\rg,$$
where $o(g_1)=\cdots=o(g_k)=3$. Then $o(g_{k+1})=\cdots =o(g_{\ld(G)})=2$.
It is clear that $\lg g_i\rg\cap \lg g_j\rg=1$ for $i\neq j$.
Thus $$|G|=1+k(3-1)+(\ld(G)-k)=\ld(G)+k+1.$$
It follows that $\ld(G)=|G|-k-1$. \qed

\medskip
By Proposition \ref{Prop=o(g_1)=3} (1), we assume that $\ld(G)< |G|-1$ in following.

\begin{prop}\label{Prop=the bound of g1}
$o(g_1)\leq |G|-\ld(G)+1$.
\end{prop}

\demo By Lemma \ref{Lemma=ld(G)=the number of maximal cyclic sub}, $\{\lg g_i\rg\mid i=1, \cdots, \ld(G)\}$
is an irredundant covering of $G$.
We have $g_2\cdots, g_{\ld(G)}$ are different pairwise elements of $G\backslash\lg g_1\rg$.
Hence $\ld(G)-1\leq |G|-o(g_1)$. It follows that $o(g_1)\leq |G|-\ld(G)+1$.
 \qed

\begin{thm}\label{Theorem=o(g_1)=|G|-ld(G)+1}
Let $G$ be a finite non-cyclic group. If $o(g_1)=|G|-\ld(G)+1$, then \\
{\rm (1)} $o(g_2)=\cdots=o(g_{\ld(G)})=2$.\\
{\rm (2)} $G\cong \D_{2n}$,
where $n=|G|-\ld(G)+1$.
\end{thm}

\demo (1) Since $G$ is non-cyclic, $\ld(G)\geq 2$.
If $o(g_2)>2$, then, by $o(g_1)=|G|-\ld(G)+1$, $$|\lg g_1\rg\cup\lg g_2\rg|\geq |G|-\ld(G)+1+2=|G|-\ld(G)+3.$$
Since $g_3, \cdots, g_{\ld(G)}$ are $\ld(G)-2$ elements in $G\backslash (\lg g_1\rg\cup\lg g_2\rg)$,
$$\ld(G)-2\leq |G|-|\lg g_1\rg\cup\lg g_2\rg|\leq |G|-(|G|-\ld(G)+3)=\ld(G)-3,$$
This is a contradiction.
Hence $o(g_2)=2$.
Thus, $o(g_3)=\cdots=o(g_{\ld(G)})=2$.

(2) Since $\ld(G)<|G|-1$, $n>2$.
It follows by (1) that $\lg g_1\rg$ is the unique cyclic subgroup of order $n$.
Take any $x\in G\backslash \lg g_1\rg$.
Then $o(x)=2$. Since $g_1x\not\in\lg g_1\rg$, $o(g_1x)=2$,
i.e. $g_1xg_1x=1$.
It follows that $g_1^x=g_{1}^{-1}$ and $\lg g_1, x\rg\cong \D_{2n}$.
We assert that $G=\lg g_1, x\rg$.
If not, then there exists an element $y\in G\backslash \lg g_1, x\rg$.
Since $y\not\in\lg g_1\rg$, $o(y)=2$. Thus $g_1^y=g_{1}^{-1}$.
It follows that $g_{1}^{xy}=g_1$.
Thus $g_1$ and $xy$ are commute and $g_1xy$ is an element of order $n$ of $G$.
This is contrary to the uniqueness of $\lg g_1\rg$.
Hence, $G=\lg g_1, x\rg\cong \D_{2n}$.\qed

%

\begin{prop}\label{Prop=o(g_i)=p}
Let $p$ be a prime. If there exists a maximal cyclic subgroup of order $p$
in $G$, then $\ld(G)\leq |G|-k(p-2)-1$, where $k$ is the number of maximal cyclic subgroups of order $p$ of $G$.
\end{prop}

\demo Let $A=\lg g_{i_1}\rg\cup\lg g_{i_2}\rg\cup\cdots\cup\lg g_{i_k}\rg$, where $\lg g_{i_1}\rg, \lg g_{i_2}\rg, \cdots, \lg g_{i_k}\rg$
are all the maximal cyclic subgroups of order $p$ of $G$. Then $|A|=1+k(p-1)$.
Since there exist $\ld(G)-k$ elements in $G\backslash A$, $\ld(G)-k\leq |G|-(1+k(p-1))$.
It follows that $\ld(G)\leq |G|-k(p-2)-1$. \qed

\begin{thm}\label{Theorem=solvable}
Let $G$ be a finite group. If $\ld(G)=|G|-t$, where $t\leq 5$, then $G$
is solvable.
\end{thm}

\demo By Proposition \ref{Prop=the bound of g1}, $o(g_1)\leq t+1$.
Since $t\leq 5$, $o(g_1)\leq 6$. If $o(g_1)\leq 4$, then $G$ is a $\{2,3\}$-group.
It follows by $p^aq^b$-theorem that $G$ is solvable. If $o(g_1)=6$, then
$t=5$ and $o(g_2)=\cdots=o(g_{\ld(G)})=2$
by Theorem \ref{Theorem=o(g_1)=|G|-ld(G)+1}. Hence
$G$ is a $\{2,3\}$-group. Hence $G$ is also solvable.

If $o(g_1)=5$, then $t=4$ or $5$ by Proposition \ref{Prop=the bound of g1}.
We assert that $t=4$. If not, then $t=5$ and $\ld(G)=|G|-5$.
Notice that $|\lg g_1\rg\cup\lg g_2\rg|=o(g_2)+(5-1)=o(g_2)+4$.
Since $g_3, \cdots, g_{\ld(G)}\in G\backslash (\lg g_1\rg\cup\lg g_2\rg)$,
$\ld(G)-2\leq |G|-|\lg g_1\rg\cup\lg g_2\rg|$.
It follows by $\ld(G)=|G|-5$ that $o(g_2)\leq 3$.
If $o(g_2)=2$, then $o(g_3)=\cdots=o(g_{\ld(G)})=2$.
Since $G=\lg g_1\rg\cup\lg g_2\rg\cup\lg g_3\rg\cup\cdots\cup\lg g_{\ld(G)}\rg$,
$|G|=5+\ld(G)-1=\ld(G)+4$. This is contrary to $\ld(G)=|G|-5$. Hence, $o(g_2)=3$.
It follows that $G$ has only one cyclic subgroup $\lg g_1\rg$ of order $5$.
Hence, $\lg g_1\rg\unlhd G$.
Clearly, $\lg g_1\rg\lg g_2\rg$ is a subgroup of order $15$ of $G$.
Notice that the group of order $15$ is cyclic. However, the largest order of maximal cyclic subgroups of $G$ is 5.
This is a contradiction. Hence, $t=4$, and hence $\ld(G)=|G|-4$.
By Proposition \ref{Prop=the bound of g1},
$o(g_2)=\cdots=o(g_{\ld(G)})=2$. Hence $|G|=2^k5$.
It follows that $G$ is solvable. \qed

\begin{rem}
In fact, if $t\leq 4$, then $G$ is supersolvable by Theorem \ref{thm=lambda(G)=|G|-2}, \ref{thm=lambda(G)=|G|-3}, \ref{thm=lambda(G)=|G|-4} in next section.
However, if $t=5$, then $G$ is not necessarily supersolvable by Theorem \ref{thm=lambda(G)=|G|-5}.
\end{rem}

\section{Finite groups with $\lambda(G)=|G|-t$ $(t\leq 5)$}

\ \ \ In this section, we give a classification of finite groups with $\lambda(G)=|G|-t(t\leq 5)$.
In following, we always assume that for $i=1,2,\cdots, {\lambda(G)}$, $g_i$ is as the same as the assume in the begining of Section 3.

\begin{thm}\label{thm=lambda(G)=|G|-2}
Let $G$ be a finite group. Then

{\rm (1)} $\ld(G)=|G|-1$ if and only if $G\cong \C_2^n$.

{\rm (2)} $\lambda(G)=|G|-2$ if and only if $G\cong \C_3$ or $S_3$.
\end{thm}

\demo (1) It is straightforward.

(2) $(\Longleftarrow).$  It is straightforward by a simple calculation. 

$(\Longrightarrow).$ By Proposition \ref{Prop=the bound of g1}, we have $o(g_1)\leq 3$.
If $o(g_1)=2$, then $G\cong \C_2^n$. It is clear that $\ld(G)=|G|-1$,
which is contrary to $\ld(G)=|G|-2$. Hence, $o(g_1)=3$.
If $G$ is cyclic, then $G\cong \C_3$. Assume $G$ is non-cyclic. Since
$o(g_1)=|G|-\ld(G)+1$, $G\cong D_6\cong S_3$ by
Theorem \ref{Theorem=o(g_1)=|G|-ld(G)+1}.\qed

\begin{thm}\label{thm=lambda(G)=|G|-3}
Let $G$ be a finite group. Then $\lambda(G)=|G|-3$ if and only if $G\cong \C_4$ or $\D_8$.
\end{thm}

\demo $(\Longleftarrow).$ It is straightforward  by a simple calculation. 

$(\Longrightarrow).$
If $G$ is cyclic, then $G\cong \C_4$. Assume that $G$ is non-cyclic.
By Proposition \ref{Prop=the bound of g1}, we have $o(g_1)\leq 4$.
We assert that $o(g_1)=4$.
If not, then $o(g_1)=2$ or $3$.
If $o(g_1)=2$, then $G\cong \C_2^n$. It is clear that $\ld(G)=|G|-1$.
This contradicts $\ld(G)=|G|-3$.
If $o(g_1)=3$, then, by Proposition \ref{Prop=o(g_1)=3},
$\ld(G)=|G|-k-1$, where $k$ is the number of subgroups of $G$ of order $3$.
Since $\ld(G)=|G|-3$, $k=2$. By Sylow theorem, $k\equiv 1(\mod 3)$. This contradicts $k=2$.
Therefore, $o(g_1)=4$. Since $o(g_1)=|G|-\ld(G)+1$, $G\cong \D_8$ by Theorem \ref{Theorem=o(g_1)=|G|-ld(G)+1}.\qed

\begin{thm}\label{thm=lambda(G)=|G|-4}
Let $G$ be a finite group. Then $\lambda(G)=|G|-4$ if and only if $G\cong \C_5$, $\D_{10}, \C_4\times \C_2$ or $\D_8\times \C_2$.
\end{thm}

\demo $(\Longleftarrow).$ It is clear that $\ld(\C_5)=1$. By Lemma \ref{Lemma=ld(C_p times C_p^k)ld(D_2n)ld(Q_2^n)} and Lemma \ref{Lemma=ld(D_8times C_2)}, $\ld(\C_4\times \C_2)=4$, $\ld(\D_{10})=6$ and  $\ld(\D_8\times \C_2)=12$.
Hence $\ld(G)=|G|-4$. 

$(\Longrightarrow).$ If $G$ is cyclic, then $G\cong \C_5$. Assume that $G$ is non-cyclic.
By Proposition \ref{Prop=the bound of g1}, we have $o(g_1)\leq 5$. We assert $o(g_1)=4$ or $5$. In fact,
if $o(g_1)=2$, then $G\cong \C_2^n$. It is clear that $\ld(G)=|G|-1$.
This contradicts $\ld(G)=|G|-4$.  If $o(g_1)=3$.
then, by Proposition \ref{Prop=o(g_1)=3},
$\ld(G)=|G|-k-1$, where $k$ is the number of subgroups of $G$ of order $3$.
Since $\ld(G)=|G|-4$, $k=3$. By Sylow Theorem, $k\equiv 1(\mod 3)$. This contradicts $k=3$.
Therefore, $o(g_1)=4$ or $5$.

If $o(g_1)=5$, then $o(g_1)=|G|-\ld(G)+1$. Thus $G\cong D_{10}$ by
Theorem \ref{Theorem=o(g_1)=|G|-ld(G)+1}.

If $o(g_1)=4$, then we will prove that $G\cong \C_4\times\C_2$ or $\D_8\times \C_2$
by three steps.

\smallskip
Step 1. $G$ is a $2$-group. 
\smallskip

It is enough to show that $o(g_i)\neq 3$ for $2\leq i\leq \ld(G)$.
Suppose that there exists a maximal cyclic subgroup of order $3$ in $G$. Since $o(g_1)=4$,
each subgroup of order $3$ is maximal cyclic in $G$.
By Proposition \ref{Prop=o(g_i)=p},
$\ld(G)\leq |G|-k-1$, where $k$ is the number of maximal cyclic subgroups of order $3$ of $G$.
Since $\ld(G)=|G|-4$, $k\leq 3$. By Sylow theorem, $k\equiv 1(\mod 3)$. It
follows that $k=1$, that is, $G$ has unique subgroup $\lg g_i\rg$ of order $3$.
Hence, $\lg g_i\rg\unlhd G$.

Since $o(g_1)=4$, $C_G(g_i)=\lg g_i\rg$. It follows that 
$G/C_G(g_i)=G/\lg g_i\rg\lesssim \Aut(\lg g_i\rg)\cong \C_2$.
Thus $|G|\mid 6$. This contradicts $o(g_1)=4$.
Hence there is no subgroup of order $3$ in $G$. Thus $G$
is a $2$-group.

\smallskip
Step 2. $G$ has only two maximal cyclic subgroups of order $4$ and the intersection
of them has order $2$.
\smallskip

Assume that $G$ has $k$ maximal cyclic subgroups $\lg g_1\rg, \cdots, \lg g_k\rg$ of order $4$.
Then $|\lg g_1\rg\cup\cdots\cup\lg g_k\rg|\geq 2k+2$.
Since $g_{k+1}, \cdots, g_{\ld(G)}$ are different elements in
$G\backslash (\lg g_1\rg\cup\cdots\cup\lg g_k\rg)$,
$$\ld(G)-k\leq |G|-|\lg g_1\rg\cup\cdots\cup\lg g_k\rg|\leq |G|-(2k+2).$$
By $\ld(G)=|G|-4$, we have $k\leq 2$.

We assert that $k=2$. If not, then $k=1$.
By step 1, $o(g_i)=2$ for $2\leq i\leq \ld(G)$.
Since $$G=\lg g_1\rg\cup\cdots\cup\lg g_{\ld(G)}\rg,\eqno (*)$$
it follows that $|G|=4+\ld(G)-1=\ld(G)+3$.
This is contrary to $\ld(G)=|G|-4$. Hence, $k=2$.

Assume that $\lg g_1\rg\cap\lg g_2\rg=1$.
Noticing that $$o(g_1)=o(g_2)=4, o(g_3)=\cdots=o(g_{\ld(G)})=2,$$ we have
$$|G|=4+3+\ld(G)-2=\ld(G)+5.$$
This is a contradiction. Hence $|\lg g_1\rg\cap\lg g_2\rg|=2$.

\smallskip
Step 3. Final results.
\smallskip

By Step 2, we have $\exp(G)=4$.
If $G$ is abelian, then $G\cong \C_4^m\times\C_2^n$.
Since there are $6$ cyclic subgroups of order $4$ in $\C_4\times\C_4$, $m=1$.
Hence, $G\cong \C_4\times\C_2^n$.
Let $k$ be the number of cyclic subgroups of order $4$ of $G$. Then
$$k=\frac{|\Omega_2(G)|-|\Omega_1(G)|}{2}=\frac{2^{n+2}-2^{n+1}}{2}=2^n.$$
It follows by step 2 that $n=1$. Thus $G\cong\C_4\times\C_2$.

Assume that $G$ is non-abelian.
Notice that
$$G=\lg g_1\rg\cup\lg g_2\rg\cup\cdots\cup\lg g_{\ld(G)}\rg,$$
where $o(g_1)=o(g_2)=4, g_1^2=g_2^2$ and $o(g_3)=\cdots=o(g_{\ld(G)})=2$.
We have $\mho_1(G)=\lg g_1^2\rg$. Since $G$ is a non-abelian $2$-group,
$G'=\mho_1(G)=\Phi(G)=\lg g_1^2\rg$.
By \cite[Lemma 4.2]{Ber1} we have
$$G\cong A_1*A_2*\cdots*A_sZ(G), \ \ {\mbox {where}}\ \  |A_i|=8.$$
Since there exist three cyclic subgroups of order $4$ in $\Q_8$,
$A_i\ncong \Q_8$. It follows that $A_i\cong \D_8$.

If $s\geq 2$, then there is a subgroup $H$ is isomorphic to $\D_8*\D_8$.
Clearly, the number of cyclic subgroups of order $4$ of $\D_8*\D_8$ is greater than $2$.
This is a contradiction. Hence $s=1$.
We assert that $\exp(Z(G))=2$. If not,
then there is a subgroup $H$ is isomorphic to $\D_8*\C_4$.
Clearly, the number of cyclic subgroups of order $4$ of $\D_8*\C_4$ is greater than $2$.
This is a contradiction. Hence $\exp(Z(G))=2$.
It follows that $G\cong \D_8\times\C_2^k$. By step 2, $G$ has only two cyclic subgroups of order $4$.
Hence $k=1$. Thus $G\cong \D_8\times\C_2$. \qed

\medskip
Before we classify finite groups with $\ld(G)=|G|-5$,  we need following lemmas.

\begin{lem}\label{o(g)_1=3,4,6}
If $\ld(G)=|G|-5$, then $o(g_1)=3, 4$ or $6$.
\end{lem}

\demo By Proposition \ref{Prop=the bound of g1}, $o(g_1)\leq 6$.
Since $\ld(\C_2^n)=2^n-1$, $o(g_1)\neq 2$.
We assert that $o(g_1)\neq 5$. If not, then $|\lg g_1\rg\cup\lg g_2\rg|=o(g_2)+(5-1)=o(g_2)+4$.
Since $g_3, \cdots, g_{\ld(G)}\in G\backslash (\lg g_1\rg\cup\lg g_2\rg)$,
$\ld(G)-2\leq |G|-|\lg g_1\rg\cup\lg g_2\rg|$.
It follows by $\ld(G)=|G|-5$ that $o(g_2)\leq 3$.
If $o(g_2)=2$, then $o(g_3)=\cdots=o(g_{\ld(G)})=2$.
Since $G=\lg g_1\rg\cup\lg g_2\rg\cup\lg g_3\rg\cup\cdots\cup\lg g_{\ld(G)}\rg$,
$|G|=5+\ld(G)-1=\ld(G)+4$. This contradicts $\ld(G)=|G|-5$. If $o(g_2)=3$, then
$G$ has only one cyclic subgroup $\lg g_1\rg$ of order $5$.
Hence, $\lg g_1\rg\unlhd G$.
Clearly, $\lg g_1\rg\lg g_2\rg$ is a subgroup of order $15$ of $G$.
Notice that the group of order $15$ is cyclic.
However, the largest order of maximal cyclic subgroups of $G$ is 5.
This is a contradiction.
Thus $o(g_1)=3, 4$ or $6$. \qed

\begin{lem}\label{Lemma=ld(G)=|G|-5 o(g_1)=3}
If $\lambda(G)=|G|-5$ and $o(g_1)=3$, then $G\cong A_4$, $\C_3^2$ or
$\langle a, b, c\mid a^3=b^3=c^2=1, a^c=a^{-1}, b^c=b^{-1}, a^b=a\rangle\cong \C_3^2\rtimes \C_2$.
\end{lem}

\demo By Proposition \ref{Prop=o(g_1)=3}, there are four cyclic subgroups of order $3$.
Let $P\in\Syl_3(G)$ and $|P|=3^k$. Since $o(g_1)=3$, $\exp (P)=3$.
Let $n_3$ be the number of subgroup of order $3$ of $P$.
Then $n_3=\frac{|P|-1}{3-1}=\frac{3^k-1}{2}$.
Clearly $n_3\leq 4$. It follows that $k\leq 2$.
Thus, $P\cong \C_3$ or $\C_3^2$.

\smallskip
{\bf Case 1.} $P\cong \C_3$.
\smallskip

Since $G$ has four cyclic subgroups of order $3$, $|G:N_G(P)|=4$ by Sylow theorem.
Let $H=N_G(P)$ and $H_G=\bigcap\limits_{x\in G}H^x$ be the core of $H$.
If $3\mid |H_G|$, then $P\leq H_G$. Thus, $G=N_G(P)H_G$ by Frattini argument.
Since $H=N_G(P)$, $H_G\leq N_G(P)$.
It follows that $G=N_G(P)$.
Hence, $P\unlhd G$. This is a contradiction.
Therefore, $3\nmid |H_G|$.
Since $G$ is a $\{2,3\}$-group and $3\nmid |H_G|$, $H_G$ is a $2$-subgroup of $G$.
Noticing that $H_G\leq H=N_G(P)$, we have $H_GP=H_G\times P$.
It follows that $H_G=1$ by $o(g_1)=3$.
Now, consider the transform representation of $H$. Since $|G:H|=4$ and $H_G=1$,
$G/H_G=G\lesssim S_4$. Noticing that $3\mid |G|$ and $4\mid |G|$, it
follows that $G\cong A_4$ or $S_4$.
By $o(g_1)=3$, we have $G\ncong S_4$. Therefore, $G\cong A_4$.

\smallskip
{\bf Case 2.} $P\cong \C_3^2$.
\smallskip

Since $o(g_1)=3$, $G$ is a $\{2, 3\}$-group.
If $2\nmid |G|$, then $G\cong \C_3^2$. Assume that $2\mid |G|$.
By Proposition \ref{Prop=o(g_1)=3}, $G$ has only four cyclic subgroups of order $3$.
It follows that $P\unlhd G$. Notice that $o(g_1)=3$.
It follows that $G/P$ is an element abelian $2$-group.
Clearly,
there is no $3'$-elements in $C_G(P)$. Thus $C_G(P)=P$.
By the $N/C$-theorem, $G/P$ is isomorphic to an elementary abelian $2$-subgroup of $GL(2,3)$.
Let
$$\a=\left(
\begin{array}{cc}
1&1\\
1&0
\end{array}\right),
\b=\left(
\begin{array}{cc}
1&-1\\
0&-1
\end{array}\right),$$
By computation, we get $o(\a)=8, o(\b)=2$ and $\a^{\b}=\a^3$. Let
$$S=\lg \a, \b\mid \a^8=\b^2=1, \a^{\b}=\a^3\rg.$$ Then $S\cong \SD_{16}.$
Hence $S$ is a Sylow $2$-subgroup of $GL(2,3)$.
By computation, $\a^4, \b, \a^2\b, \a^4\b, \a^6\b$ are all the elements of order $2$ of $S$.

Let $P=\lg a,b\rg$ and $1\neq c\in Q$, where $Q\in{\Syl}_2(G)$. Since the largest order of elements of $G$ is 3, $c$ is an element of order $2$ of $G$.
Since $P\unlhd G$, $c$ induces an automorphism of $P$ by conjugate.
Let $a^c=a^ib^j$ and $b^c=a^kb^l$. Then the map
$$\sigma: c\mapsto \left(\begin{array}{cc}
            i & j \\
            k & l
          \end{array}\right)$$
is an isomorphic mapping from $Q$ to an elementary abelian $2$-subgroup of $S$.
Since $o(g_1)=3$, $C_P(c)=1$ for any element $1\neq c\in Q$.
By computation, $\sigma(Q)=\langle \alpha^4\rangle$.
Since $\a^4=\left(\begin{array}{cc}
            -1 & 0 \\
            0 & -1
          \end{array}\right)$,
$a^c=a^{-1}$ and $b^c=b^{-1}$.
It follows that $$G=\langle a, b, c\mid a^3=b^3=c^2=1, a^c=a^{-1}, b^c=b^{-1}, a^b=a\rangle.$$
\qed

%
%
%

\begin{lem}\label{Lemma=ld(G)=|G|-5 o(g_1)=4}
If $\lambda(G)=|G|-5$ and $o(g_1)=4$, then $G\cong \Q_8$.
\end{lem}

\demo We prove the result by the following four steps.

\smallskip
Step 1. $G$ is a $2$-group.
\smallskip

If not, then there exist maximal cyclic subgroups of order $3$ in $G$.
Assume that $G$ has $k$ maximal cyclic subgroups of order $3$.
Let $A$ be the union of the maximal cyclic subgroups of order $3$ of $G$.
Then $|A|=1+2k$. Since there are at least $\ld(G)-(k+1)$ elements in $G\backslash (\lg g_1\rg\cup A)$,
we have $$\ld(G)-(k+1)\leq |G|-|\lg g_1\rg\cup A|=|G|-(2k+4).$$
It follows by $\ld(G)=|G|-5$ that $k\leq 2$. By Sylow theorem, we have $k=1$.

Let $H$ be the maximal cyclic subgroup of order $3$ of $G$. Then $H\unlhd G$.
Since $o(g_1)=4$, $C_G(H)=H$.
It follows that $G/H\lesssim \Aut(H)\cong\C_2$. Thus $|G|=3$ or $6$.
This is contrary to $o(g_1)=4$.
Therefore, $G$ has no subgroup of order $3$.
Since $o(g_1)=4$, $G$ is a $2$-group.

\smallskip
Step 2.  $G$ has only three maximal cyclic subgroups of order $4$.
\smallskip

Let $\lg g_1\rg, \lg g_2\rg, \cdots, \lg g_k\rg$ be the maximal cyclic subgroups of order $4$ of $G$.
Then $|\lg g_1\rg\cup\lg g_2\rg\cup \cdots\cup \lg g_k\rg|\geq 2k+2$.
Since $g_{k+1}, \cdots, g_{\ld(G)}$ are the elements of $G\backslash (\lg g_1\rg\cup\lg g_2\rg\cup \cdots\cup \lg g_k\rg)$,
$$\ld(G)-k\leq |G|-|\lg g_1\rg\cup\lg g_2\rg\cup \cdots\cup \lg g_k\rg|\leq |G|-(2k+2).$$
It follows by $\ld(G)=|G|-5$ that $k\leq 3$.

If $k=1$, then $o(g_2)=\cdots=o(g_{\ld(G)})=2$.
By $G=\lg g_1\rg\cup\lg g_2\rg\cup\cdots\cup\lg g_{\ld(G)}\rg$,
we have $|G|=4+\ld(G)-1=\ld(G)+3$, which is contrary to $\ld(G)=|G|-5$.
If $k=2$, then $o(g_1)=o(g_2)=4$ and $o(g_3)=\cdots=o(g_{\ld(G)})=2$.
By $G=\lg g_1\rg\cup\lg g_2\rg\cup\cdots\cup\lg g_{\ld(G)}\rg$,
we have $|G|=|\lg g_1\rg\cup\lg g_2\rg|+\ld(G)-2$.
Since $\ld(G)=|G|-5$, $|\lg g_1\rg\cup\lg g_2\rg|=7$.
It follows that $\lg g_1\rg\cap\lg g_2\rg=1$.
By step 1, $G$ is a $2$-group. Hence $Z(G)\neq 1$.
Take an element $a$ of order $2$ in $Z(G)$. Then $a\not\in \lg g_1\rg$
or $a\not\in \lg g_2\rg$.
Without losing generality, assume that $a\not\in \lg g_1\rg$.
 Then $\lg ag_1\rg$ is a cyclic subgroup of order $4$.
 Since $k=2$ and
 $a\not\in \lg g_1\rg$, $\lg ag_1\rg=\lg g_2\rg$.
 Thus $g_2^2=(ag_1)^2=g_1^2\in\lg g_1\rg\cap\lg g_2\rg$,
 which is contrary to $\lg g_1\rg\cap\lg g_2\rg=1$.
 To sum up, $k=3$. That is, $G$ has only three maximal cyclic subgroups $\lg g_1\rg,
\lg g_2\rg$ and $\lg g_3\rg$ of order $4$.

\smallskip
Step 3.  $\lg g_1,g_2,g_3\rg\cong Q_8$.
\smallskip

If $|\lg g_1\rg\cap\lg g_2\rg\cap\lg g_3\rg|=1$,
then $|\lg g_1\rg\cup\lg g_2\rg\cup\lg g_3\rg|\geq 9$.
Noticing that $G=\lg g_1\rg\cup\lg g_2\rg\cup\cdots\cup\lg g_{\ld(G)}\rg$
and $o(g_4)=\cdots=o(g_{\ld(G)})=2$,
we have $$|G|=|\lg g_1\rg\cup\lg g_2\rg\cup\lg g_3\rg|+\ld(G)-3\geq 9+\ld(G)-3=\ld(G)+6.$$
It follows that $\ld(G)\leq |G|-6$, which is contrary to $\ld(G)=|G|-5$.
Hence $|\lg g_1\rg\cap\lg g_2\rg\cap\lg g_3\rg|\neq 1$
and $g_1^2=g_2^2=g_3^2$.

We claim that $\lg g_1,g_2,g_3\rg$ is not abelian. Otherwise, $\lg g_1g_2g_3\rg$ is the 4th maximal cyclic subgroup of order 4, a contradiction.  Without loss of generality, we assume that $[g_1,g_2]\ne 1$. It is easy to see that $\lg g_1,g_2\rg\cong Q_8$ and $g_3=g_1g_2$.

\smallskip
Step 4. $G\cong Q_8$.
\smallskip

By $$G=\lg g_1\rg\cup\lg g_2\rg\cup\cdots\cup\lg g_{\ld(G)}\rg$$
and
$$o(g_1)=o(g_2)=o(g_3)=4, o(g_4)=\cdots=o(g_{\ld(G)})=2, g_1^2=g_2^2=g_3^2 ,$$
we have $\Phi(G)=\mho_1(G)=\lg g_1^2\rg$.
It follows that $\lg g_1\rg\unlhd G$.

We assert $C_G(g_1)=\lg g_1\rg$.
If not, there exists an element $a\in C_G(g_1)\backslash \lg g_1\rg$.
It is easy to see that $\lg g_1a\rg $ is the 4th maximal cyclic subgroup of order 4, a contradiction.
Hence $C_G(g_1)=\lg g_1\rg$.
By $N/C$-theorem, $G/C_G(g_1)=G/\lg g_1\rg\lesssim C_2$.
Thus, $|G|\leq 8$. By Step 3, $G=\lg g_1,g_2,g_3\rg\cong \Q_8$.
\qed

\begin{lem}\label{Lemma=ld(G)=|G|-5 o(g_1)=6}
If $\lambda(G)=|G|-5$ and $o(g_1)=6$, then $G\cong\C_6$ or $\D_{12}$.
\end{lem}

\demo If $G$ is cyclic, then $\ld(G)=1$.
Since $\ld(G)=|G|-5$, we get $G\cong\C_6$. Assume that $G$ is non-cyclic.
Since
$o(g_1)=|G|-\ld(G)+1$, $G\cong \D_{12}$ by
Theorem \ref{Theorem=o(g_1)=|G|-ld(G)+1}. \qed

\begin{thm}\label{thm=lambda(G)=|G|-5}
Let $G$ be a finite group. Then $\lambda(G)=|G|-5$ if and only if $G$ is isomorphic to one of the groups in Lemma \ref{Lemma=ld(G)=|G|-5 o(g_1)=3}, \ref{Lemma=ld(G)=|G|-5 o(g_1)=4} and \ref{Lemma=ld(G)=|G|-5 o(g_1)=6}.
\end{thm}

\demo $(\Longrightarrow).$ It follows by Lemma \ref{o(g)_1=3,4,6}--\ref{Lemma=ld(G)=|G|-5 o(g_1)=6}.
%

$(\Longleftarrow).$ It is clear that $\ld(\C_6)=1$.
By Lemma \ref{Lemma=ld(C_p times C_p^k)ld(D_2n)ld(Q_2^n)},
$\ld(\Q_8)=3$, $\ld(\D_{12})=7$ and  $\ld(\C_3^2)=4$.
It is clear that the order of maximal cyclic subgroup of $A_4$ is $2$ or $3$.
Since there are three cyclic subgroups of order $2$ and four cyclic subgroups of order $3$,
$\ld(A_4)=7$.

Assume that $G=\langle a, b, c\mid a^3=b^3=c^2=1, a^c=a^{-1}, b^c=b^{-1}, a^b=a\rangle$.
Let $P=\lg a, b\rg$. Then $P\in\Syl_3(G)$ and $P\unlhd G$.
It is clear that $x^c=x^{-1}, \forall x\in P$. Thus $N_G(\lg c\rg)=\lg c\rg$
and the order of maximal cyclic subgroups of $G$ is $2$ or $3$.
It follows that $G$ has four cyclic subgroups of order $3$ and nine cyclic subgroups of order $2$.
Thus $\ld(G)=13$.

To sum up,  $\ld(G)=|G|-5$ for the groups listed in Lemma \ref{Lemma=ld(G)=|G|-5 o(g_1)=3},
\ref{Lemma=ld(G)=|G|-5 o(g_1)=4} and \ref{Lemma=ld(G)=|G|-5 o(g_1)=6}. \qed

\end{document}